\documentclass[11pt]{amsart}
\usepackage{amssymb,pstricks,pst-plot}

\definecolor{OrangeRed}{cmyk}{0,0.6,1,0}
\definecolor{DarkBlue}{cmyk}{1,1,0,0.20}
\definecolor{Black}{cmyk}{0,0,0,1}
\definecolor{Violet}{cmyk}{0.79,0.88,0,0}
\definecolor{myblue}{rgb}{0.66,0.78,1.00}

\parskip=\smallskipamount

\newtheorem{theorem}{Theorem}[section]
\newtheorem{lemma}[theorem]{Lemma}
\newtheorem{corollary}[theorem]{Corollary}
\newtheorem{proposition}[theorem]{Proposition}

\theoremstyle{definition}
\newtheorem{definition}[theorem]{Definition}

\newtheorem{remark}[theorem]{Remark}

\newcommand{\C}{\mathbb{C}}

\newcommand{\N}{\mathbb{N}}

\newcommand{\Z}{\mathbb{Z}}

\newcommand{\R}{\mathbb{R}}

\newcommand{\cJ}{\mathcal{J}}

\newcommand{\cO}{\mathcal{O}}

\newcommand\wt{\widetilde}
\newcommand\spsh{strongly plurisubharmonic}

\def\di{\partial}

\def\bs{\backslash}

\renewcommand\i{{\mathrm{i}}}

\newcommand\CAP{{\rm CAP}}
\newcommand\PCAP{{\rm PCAP}}

\numberwithin{equation}{section}

\begin{document}
\title[The Oka principle for sections of  stratified fiber bundles]
{The Oka principle for sections of \\ stratified fiber bundles}
\author{Franc Forstneri\v c}
\address{Faculty of Mathematics and Physics, University of Ljubljana, 
and Institute of Mathematics, Physics and Mechanics, Jadranska 19, 
1000 Ljubljana, Slovenia}
\email{franc.forstneric@fmf.uni-lj.si}
\thanks{Supported by the research program P1-0291, Republic of Slovenia.}

%
%
\date{March 14, 2009}

\begin{abstract}
Let $Y$ be a complex manifold with the property that
every holomorphic map from a neighborhood of a  compact convex
set $K$ in a complex Euclidean space $\C^n$ to $Y$ 
can be approximated, uniformly on $K$,
by entire maps $\C^n\to Y$. If $X$ is a reduced Stein space
and $\pi \colon Z\to X$ is a holomorphic fiber bundle
with fiber $Y$ then we show that sections $X\to Z$ satisfy the
Oka principle with approximation and interpolation.
The analogous result is proved for sections of stratified fiber bundles,
and of submersions with stratified sprays over Stein spaces.
\end{abstract}
\maketitle

\centerline{\em Dedicated to Joseph J.\ Kohn on the occasion
of his 75th birthday}

\section{Introduction}
In this paper we establish the {\em Oka principle} for sections 
of stratified holomorphic fiber bundles over Stein spaces
under the condition that all fibers satisfy a Runge type
approximation property introduced in \cite{FF:CAP,FF:EOP},
and also for sections of submersions with stratified sprays.
Our main results are Theorems \ref{CAP:stratified}, 
\ref{Oka:parametric} and \ref{stratified-sprays}.
All complex spaces will be assumed reduced and paracompact.
For the general theory of Stein spaces we 
refer to \cite{Grauert-Remmert0}.

We recall from \cite{FF:CAP} the definition of the
{\em Convex Approximation Property} of a complex manifold;
this will be the main assumption on the fibers of the 
stratified fiber bundles that we shall consider.
Let $z=(z_1,\ldots,z_n)$ be complex coordinates on $\C^n$, with
$z_j=x_j+\i y_j$. Set
\begin{equation}
\label{cube}
    Q =\{(z_1,\ldots,z_n)\in \C^n\colon |x_j| \le 1,\ |y_j|\le 1,\ j=1,\ldots,n\}.
\end{equation}
A {\em special convex set} in $\C^n$ is a compact convex set $K$ of the form
\begin{equation}
\label{special-convex}
    K =\{z\in Q \colon y_n \le h(z_1,\ldots,z_{n-1},x_n)\},
\end{equation}
where $Q$ is the cube (\ref{cube}) and
$h$ is a continuous concave function with values
in $(-1,1)$. Such $(K,Q)$ will be called a
{\em special convex pair} in $\C^n$.

Given a closed set $K$ in a complex space $X$, 
a map $f\colon K\to Z$ to a complex space $Z$ is said to be
holomorphic if $f$ is the restriction to $K$ of a holomorphic map 
$U\to Z$ from an open neighborhood $U$ of $K$ in $X$.
For a homotopy of maps the neighborhood $U$ 
will not depend on the parameter.

\begin{definition}  
\label{CAP}
A complex manifold $Y$ satisfies the {\em Convex Approximation Property}
{\rm (CAP)} if for every special convex pair $(K,Q)$ in $\C^n$
$(n\in\N)$ and for every holomorphic map $f\colon K\to Y$ there is 
a sequence of holomorphic maps $f_k\colon Q\to Y$ $(k=1,2,\ldots)$ 
converging to $f$ uniformly on $K$. 
\end{definition}

By Theorem \ref{CAP:stratified} below, CAP is equivalent to the condition 
that every holomorphic map $K\to Y$ from a compact convex set 
$K\subset \C^n$ $(n\in\N)$ can be approximated,
uniformly on $K$, by entire maps $\C^n\to Y$. 

Next we recall the notions of a holomorphic submersion and 
of a stratified holomorphic fiber bundle over a complex space.

%
%
\begin{definition}
\label{submersion}
Let $X$ and $Z$ be reduced complex spaces.
A holomorphic map $\pi \colon Z\to X$ is a {\em holomorphic submersion} 
if for every point $z_0\in Z$ 
there exist an open neighborhood $V\subset Z$ of $z_0$,
an open neighborhood $U\subset X$ of $x_0=\pi(z_0)\in X$,
an open set  $W$ in $\C^p$, and a biholomorphic map
$\phi\colon V\to U\times W$ such that $pr_1\circ \phi= \pi$,
where $pr_1 \colon U\times W\to U$ is the projection.
\end{definition}

A {\em stratification} of a finite dimensional complex space
$X$ is a finite descending chain of closed complex subvarieties 
\begin{equation}
\label{stratification}
    X=X_0\supset X_1\supset\cdots\supset  X_m=\emptyset
\end{equation}
such that each connected component
$S$ (stratum) of any difference $X_k\bs X_{k+1}$ is a complex 
manifold satisfying $\overline S\bs S\subset X_{k+1}$ whose dimension
$\dim S$ only depends on $k$.

\begin{definition}
\label{SFB}
A holomorphic submersion $\pi \colon Z\to X$ (Def.\ \ref{submersion}) 
is a {\em stratified holomorphic fiber bundle}
if $X$ admits a stratification (\ref{stratification})
such that the restriction  $Z|_S\to S$ to every stratum $S\subset X_k\bs X_{k+1}$
is a holomorphic fiber bundle over $S$. 
\end{definition}

Given a complex space $(X,\cO_X)$, we denote by 
$\cO(X)$ the algebra of all holomorphic functions
on $X$. Recall that a compact set $K$ in $X$ is $\cO(X)$-convex
if for every point $p\in X\bs K$ there exists
an $f\in\cO(X)$ such that $|f(p)| > \sup_{x\in K} |f(x)|$. 

We are now ready to state our first main result
which will also serve as the definition of various Oka-type 
properties; for the proof see \S \ref{S5}.

\begin{theorem}
\label{CAP:stratified}
Assume that $X$ is a reduced Stein space, that $\pi\colon Z\to X$ is a
holomorphic submersion, and that 
$X$ is exhausted by a sequence of open relatively compact 
subsets $U_1\subset U_2 \subset\cdots \subset \bigcup_{j=1}^\infty U_j=X$ such that
each restriction $Z|_{U_j}\to U_j$ is a stratified holomorphic fiber
bundle all of whose fibers enjoy \CAP (Def.\ \ref{CAP}). 
Choose a distance function $d$ on $Z$.
Then the following hold:

{\em (A) The basic Oka principle:} Every continuous section $f\colon X\to Z$ 
of $\pi\colon Z\to X$ is homotopic to a holomorphic section. 

{\em (B) The basic Oka principle with interpolation and approximation:}
Given a continuous section $f\colon X\to Z$ that is holomorphic
in a neighborhood of a compact $\cO(X)$-convex subset
$K$ in $X$ and whose restriction $f|_{X'}\colon X'\to Z$ 
to a closed complex subvariety 
$X'$ of $X$ is holomorphic on $X'$, there exists for every $\epsilon>0$ 
a homotopy of continuous sections $f_t\colon X\to Z$ $(t\in[0,1])$,
with $f_0=f$, satisfying the following properties:
\begin{itemize}
\item[(i)]     $f_1$ is holomorphic on $X$,
\item[(ii)]    $f_t|_{X'}=f_0|_{X'}$ for each $t\in [0,1]$, and
\item[(iii)]   $f_t$ is holomorphic on $K$
and $\sup_{x\in K} d\bigl(f_t(x),f_0(x)\bigr) <\epsilon$ for each $t\in [0,1]$.
\end{itemize}

{\em (C) The basic Oka principle with approximation and jet interpolation:}
Given $K$ and $X'$ as in {\rm (B)} and a continuous section 
$f\colon X\to Z$ that is holomorphic in an open set $V\supset K\cup X'$,
there exists for every $\ell \in \N$ and $\epsilon>0$ 
a homotopy $\{f_t\}_{t\in[0,1]}$ satisfying {\rm (B)} 
such that $f_t$ is holomorphic near $X'$, and it agrees with $f=f_0$ to order 
$\ell$ along $X'$ for every $t\in[0,1]$.

{\em (D) The one-parametric Oka principle:} 
Every homotopy $f_t\colon X\to Z$ $(t\in[0,1])$ of sections  between a pair of 
holomorphic sections $f_0,f_1$ can be deformed, with fixed ends at $t=0,1$,
to a homotopy consisting of holomorphic sections.
If the homotopy $f_t$ is fixed on a closed complex subvariety
$X'$ of $X$ then the deformation can be chosen fixed on $X'$. 
\end{theorem}

Combining parts (A) and (D) in the conclusion of Theorem \ref{CAP:stratified} 
we obtain the following corollary.
(See also Corollary \ref{whe} below.)

\begin{corollary}
\label{CAP-1parametric}
If $Z\to X$ is a stratified holomorphic fiber bundle over a Stein space $X$
such that all its fibers enjoy \CAP, then the inclusion 
\[
	 \iota\colon \Gamma_{\cO}(X,Z) \hookrightarrow \Gamma(X,Z)
\]
of the space of  holomorphic sections into the space of 
continuous sections induces a bijection of the path connected components
of the two spaces (endowed with the compact-open topology).
\end{corollary}

Many examples and sufficient conditions for the validity
of CAP can be found in the papers \cite{FF:CAP,FF:EOP,FF:flex}.
In particular, if $Y$ admits a dominating holomorphic spray
in the sense of Gromov \cite{Gromov} (such $Y$ is said to be {\em elliptic}), 
or, more generally, if it admits a finite dominating collection of holomorphic sprays 
(see \cite{FF:subelliptic}; such $Y$ is said to be {\em subelliptic}),
then $Y$ satisfies CAP, and also its fully parametric version, PCAP
(see Def.\ \ref{def:PCAP} below). The proof of this result is essentially a reduction
to the Oka-Weil approximation theorem; see Theorem 4.1 in \cite[p.\ 135]{FP1}
and Theorem 3.1 in \cite[p.\ 534]{FF:subelliptic}.
Hence Theorem \ref{CAP:stratified} implies the following corollary.

\begin{corollary}
\label{spray:stratified}
If $\pi\colon Z\to X$ is a stratified holomorphic fiber bundle 
over a Stein space $X$ such that all its fibers are elliptic 
in the sense of Gromov \cite{Gromov}, or subelliptic in 
the sense of \cite{FF:subelliptic}, then sections  $X\to Z$ 
satisfy the conclusions of Theorem \ref{CAP:stratified}.
\end{corollary}

We conclude this introduction by a brief survey of the known results
on the Oka principle for section of holomorphic submersions.

The classical {\em Oka-Grauert principle} is essentially 
Theorem \ref{CAP:stratified} in the  special case when 
$\pi\colon Z\to X$ is a holomorphic principal bundle with a 
Lie group fiber, or an associated holomorphic fiber bundle 
with a complex homogeneous fiber. This was proved 
by Oka \cite{Oka} for the Lie group $\C^*$, and by Grauert 
in general \cite{Grauert3}. For expositions and extensions 
see Cartan \cite{Cartan}, Grauert and Kerner \cite{GK}, 
Forster and Ramspott \cite{FRam1,FRam2}, Henkin and Leiterer \cite{HL:Oka},
and Heinzner and Kutzschebauch \cite{HK}.
Every complex homogeneous manifold  satisfies CAP
(see Grauert \cite{Grauert1, Grauert2}). 
Even for bundles with homogeneous fibers, Theorem \ref{CAP:stratified} is 
stronger than Grauert's theorem since no conditions are imposed on the
transition maps.

The case of Theorem  \ref{CAP:stratified} 
when $X$ is a Stein manifold (without singularities) and
$\pi\colon Z\to X$ is a holomorphic fiber bundle whose fiber $Y$ admits 
a dominating holomorphic spray is due to Gromov \cite{Gromov}. 
For an exposition and extensions of Gromov's work see the papers
\cite{FF:subelliptic,FF:multivalued,FP1,FP2,FP3,Prezelj3};
for a homotopy theoretic point of view see also 
L\'arusson \cite{Larusson1,Larusson2,Larusson3}.

The CAP property was first introduced in \cite{FF:CAP}
where Theorem \ref{CAP:stratified} was proved in the case when $X$ is a Stein manifold 
and $X'=\emptyset$, i.e., without the interpolation conditions (B-ii)
and (C). Interpolation  was added 
in \cite[Theorem 1.1]{FF:EOP} (see also \cite{Larusson3}). 
By \cite[Corollary 1.3]{FF:EOP} (see also L\'arusson \cite{Larusson3})
the CAP property of a complex manifold $Y$ is equivalent to several 
ostensibly different Oka properties for maps of Stein manifolds to $Y$, 
expressed by conditions (B-ii), (B-iii) and (C)
in Theorem \ref{CAP:stratified}. 
The one-parametric Oka principle (part (D) in Theorem \ref{CAP:stratified}) 
is a simple consequence of the Oka principle with interpolation, 
(B-ii). (See the proof of Theorem \ref{CAP:stratified} in \S \ref{S5}.)

Although it was remarked in \cite{FF:CAP, FF:EOP}
that Theorem \ref{CAP:stratified} also holds when
$X$ is a Stein space with singularities,
a complete proof has not been available up to now.
The outline proposed in \cite[Remark 6.6, p.\ 705]{FF:CAP}
requires certain not entirely obvious technical improvements
in order to obtain the approximation statement. 
Here we give a complete proof,  thereby dispelling the 
impressioon that this is only a theory for manifolds.

We mention that the Oka principle furnished by Theorem \ref{CAP:stratified} 
and Corollary \ref{spray:stratified}
was used in the proof of embedding theorems for Stein spaces into
Euclidean spaces of minimal dimension
(see \cite{EG,Schurmann,Prezelj2}).  
For embeddings with interpolation 
on discrete sequences see also \cite{FIKP, Prezelj1}.

In \S\ref{parametric} we consider the parametric Oka principle.
In \S\ref{S7} we obtain existence results  for holomorphic sections 
under suitable connectivity assumptions on the fibers.

Combining the induction scheme in this paper with the methods 
from \cite{FP2,FP3} we also obtain the Oka principle for sections of 
{\em submersions with stratified sprays} over Stein spaces
(see Theorem \ref{stratified-sprays} below). 
This result, which generalizes the work of Forster and Ramspott \cite{FRam2},
was first stated by Gromov \cite{Gromov}. Complete details in the nonstratified
case (i.e., for submersions with fiber dominating sprays over Stein manifolds) 
were given in \cite{FP2}. 
A sketch of proof for the stratified case can be found in \cite[\S 7]{FP3},
but without the approximation condition (only interpolation).
The additional details, given in \S 2--\S 5 of this paper, enable 
us to complete the proof of the full Oka principle as indicated 
in \cite{FP3}. This general version of the Oka principle has recently 
been used by Ivarsson and Kutzschebauch \cite{IK} in a solution 
of the {\em holomorphic Vaserstein problem} posed by 
Gromov in \cite{Gromov}.

Some recent extensions of the Oka principle should be 
mentioned: For sections of subelliptic holomorphic submersions
over 1-convex manifolds (see Prezelj \cite{Prezelj3});
for sections of holomorphic Banach bundles over 1-convex manifolds
(Leiterer and V\^aj\^aitu \cite{Leiterer-Vajaitu}); and 
the {\em soft Oka principle} to the effect
that the Oka principle holds universally if one allows 
homotopic deformations of the Stein structure on the source manifold
(see \cite{FSlapar1,FSlapar2}).

In the course of proving Theorem \ref{CAP:stratified} (see \S\ref{S5}) 
we obtain some results of  independent interest.
In \S\ref{S2} we improve \cite[Theorem 2.1]{FF:EOP}
concerning the existence of open Stein neighborhoods of the
set $K\cup X'$ in Theorem \ref{CAP:stratified}.
In \S\ref{S3} we obtain a parametric version
of a theorem of Docquier and Grauert \cite{DG}
on the existence of holomorphic retractions onto
Stein submanifolds. In \S\ref{S4} we prove
a semiglobal approximation/extension theorem
that allows passage from a stratum to the next higher one.

%
%
%
%
\section{Open Stein neighborhoods}
\label{S2}
The following result extends Theorem 2.1 from \cite{FF:EOP}.

%
%
\begin{theorem}
\label{Stein-nbds}
Let $(X,\cO_X)$ be a paracompact complex space, possibly nonreduced,
and let $X'$ be a closed Stein subvariety of $X$.
Assume that $K$ is a compact set in $X$ that is
$\cO(\Omega)$-convex in an open Stein domain $\Omega\subset X$
containing $K$ and such that $K\cap X'$ is $\cO(X')$-convex.
For every open set $U$ in $X$ containing $K\cup X'$
there exists an open Stein domain $V$, 
with $K\cup X' \subset  V\subset U$, 
such that $K$ is $\cO(V)$-convex.
\end{theorem}

The special case of Theorem \ref {Stein-nbds} with $K=\emptyset$ 
is due to Siu \cite{Siu} (see also Col\c toiu \cite{Co} and Demailly \cite{De}).
The small improvement over \cite[Theorem 2.1]{FF:EOP} is that
Stein neighborhoods $V$ of $K\cup X'$ can be chosen such that $K$ is
$\cO(V)$-convex. This property plays a key role in
holomorphic approximations theorems, and hence it makes the result
much more useful. Theorem \ref{Stein-nbds} is used in \S \ref{S4} below 
to obtain a semiglobal approximation-interpolation
result, Proposition \ref{approximation}, for sections of holomorphic submersions.

As shown in \cite{FF:EOP},
the necessity of $\cO(X')$-convexity of  $K\cap X'$
in Theorem \ref{Stein-nbds} is seen by taking 
\[
	X=\C^2, \quad X'=\C\times\{0\}, \quad 
	K=\{(z,w) \in\C^2 \colon 1\le |z|\le 2,\ |w|\le 1\}.
\]
In this example, every Stein neighborhood of $K\cup X'$
also contains the bidisc $\{(z,w)\colon |z|\le 2,\ |w|\le 1\}$.

\begin{proof}[Proof of Theorem \ref{Stein-nbds}]
A compact set $K$ that is $\cO(\Omega)$-convex
in an open Stein domain $\Omega\supset K$ in $X$ will be 
called {\em holomorphically convex}.
By the classical theory this is equivalent
to the existence of a plurisubharmonic function
$\rho_0 \colon \Omega\to\R_+$ such that $\rho^{-1}_0(0)=K$;
we can choose $\rho_0$ to be smooth and strongly plurisubharmonic on
$\Omega\bs K$. (See \cite[Theorem 5.1.5]{Ho}.)
Fix a function $\rho_0$ with these properties.

Choose an open set $U \subset X$ containing $K\cup X'$.
As $X'$ is Stein, the construction in \cite[p.\ 737]{FF:EOP} gives an open set
$W$ in $X$ with $K\cup X'\subset W \subset U$, a (small) number $c>0$,
and a smooth plurisubharmonic function $\rho\colon W\to \R_+$
that agrees with $\rho_0$ in
$U_c=\{x\in \Omega\colon \rho_0(x)<c\} \supset K$
and that satisfies $\rho>c$ on $W\bs \overline U_c$.
(See (i) and (ii) at the bottom of page 737 in \cite{FF:EOP}.)

By \cite[Theorem 2.1]{FF:EOP} there exists an open Stein domain $V$ in $X$
satisfying $K\cup X' \subset V\subset W$. (In \cite{FF:EOP},
the subvariety $X'$ was denoted $X_0$. The ambient complex space
$X$ was assumed to be reduced, but the latter property was never
used in the proof. Compare with \cite{De}.)
The restriction $\rho|_V$ is a nonnegative plurisubharmonic function
that vanishes precisely on $K$. Since $V$ is Stein,
it follows that $K$ is $\cO(V)$-convex.
\end{proof}

%
%
%
%
\section{Holomorphic retractions}
\label{S3}
The following well known result
is obtained by combining a theorem of Docquier and Grauert
(\cite[Satz 3]{DG}, \cite[p.\ 257, Theorem 8]{GR})
with Siu's theorem on the existence of  open Stein
neighborhoods of Stein subvarieties \cite[Corollary 1]{Siu}.

%
%
\begin{theorem}
\label{retraction-on-submanifold}
Let $S$ be a Stein manifold, embedded as a locally closed
complex submanifold in a complex manifold $M$.
There exist an open Stein neighborhood $\Omega$ of $S$ in $M$
and a homotopy of holomorphic maps
$\iota_t\colon \Omega \to \Omega$ $(t\in [0,1])$ satisfying
\begin{itemize}
\item[(a)]  $\iota_0$ is the identity map on $\Omega$,
\item[(b)]  $\iota_t|_S$ is the identity map on $S$ for
all $t\in [0,1]$, and
\item[(c)]  $\iota_1(\Omega)=S$.
\end{itemize}
\end{theorem}

The family $\{\iota_t\}_{t\in[0,1]}$
is a strong deformation retraction of $\Omega$ onto $S$ consisting
of holomorphic mappings. There is no immediate analogue of
Theorem \ref{retraction-on-submanifold} when $S$
is a Stein space with singularities.

We shall need the following parametric version
of Theorem \ref{retraction-on-submanifold}.

%
%

\begin{proposition}
\label{retractions-on-fibers}
Assume that $M$ and $B$ are finite dimensional complex spaces, 
$\pi\colon M\to B$ is a holomorphic submersion
(Def.\ \ref{submersion}), and $S$ is a locally closed
Stein subvariety of $M$ whose fiber $S_b= S\cap M_b$
$(b\in B)$ is a locally closed complex submanifold of $M_b$
with dimension independent of $b\in B$.
Then there exist an open Stein neighborhood $\Omega$ of $S$
in $M$ and a homotopy of holomorphic maps
$\iota_t\colon \Omega\to \Omega$ $(t\in [0,1])$ satisfying
properties (a)--(c) in Theorem \ref{retraction-on-submanifold}
and such that $\pi\circ\iota_t=\pi$ for each $t\in[0,1]$.
\end{proposition}

\begin{proof}
We begin by recalling the proof in the classical case
when $B$ is a singleton. Now $S$ is a locally closed Stein submanifold
of a complex manifold $M$. By Theorem  \ref{Stein-nbds}
we can replace $M$ by an open Stein neighborhood of $S$
that contains $S$ as a closed submanifold. We consider
the holomorphic tangent bundle $TS$ as a subbundle
of $TM|_S$, making the usual identification of $TM$
with $T^{(1,0)} M$. By Cartan's Theorem A
there exist finitely many holomorphic vector fields
$V_1,\cdots, V_N$ on $M$ such that the vectors
$V_1(x),\ldots,V_N(x)$ span $T_x M$ over $\C$ at every point
$x\in M$. Let $\phi^j_t$ denote the flow of $V_j$, i.e.,
the solution of the ordinary differential equation
$\dot \phi^j_t(x) = V_j(\phi^j_t(x))$ satisfying the initial
condition $\phi^j_0(x)=x$. For each point $x\in M$ there is a
number $T=T(x)>0$ such that $\phi^j_t(x)$ exists for all
$t\in \C$ in the disc $|t|<T$, and $T$ can be
chosen independent of the point $x$ belonging
to a compact subset of $M$. The map
\begin{equation}
\label{eqn:flows}
    F(x,t_1,\ldots,t_N)=
    \phi^1_{t_1}\circ\cdots\circ\phi^{N}_{t_N} (x),
    \quad x\in S,\ t_1,\ldots,t_N \in\C
\end{equation}
is defined and holomorphic in an open
neighborhood of $S\times \{0\}^N$ in
$S \times \C^N$. We have $F(x,0)=x$ and
$\frac{\di}{\di t_j} F(x,t)|_{t=0} = V_j(x)$
for every $x\in S$. As the $V_j$'s span $TM$ at every point,
$F$ is a submersion along $S\times \{0\}^N$. Let
\[
   \Theta_x = \di_t|_{t=0} F(x,t) \colon \C^N \to T_x M,\quad x\in S.
\]
The subset $E\subset S\times\C^N$ with the fibers
\[
    E_x = \{v\in \C^N\colon \Theta_x(v) \in T_x S \},\quad x\in S
\]
is then a holomorphic vector subbundle of the trivial bundle
$S \times \C^N$. Since $S$ is Stein, there is a holomorphic 
vector subbundle $\nu\subset S\times\C^N$ such that
$S\times\C^N = E\oplus\, \nu$ \cite[Theorem 7,\ p.\ 256]{GR}.
By the construction, $\Theta \colon \nu \to TM|_S$ is an injective
holomorphic vector bundle map and $TM|_S = TS \oplus \Theta(\nu)$;
thus $\nu \cong\ \Theta(\nu)$ is the normal
bundle of $S$ in $M$.

Consider the restriction of $F$ to $\nu$.
By the inverse function theorem $F$
maps an open neighborhood $\Omega'$ of the zero section
in $\nu$ biholomorphically onto an open neighborhood
$\Omega =F(\Omega')$ of $S$ in $M$. Further, choosing $\Omega'$
to have convex fibers, $F$ conjugates the family
of radial dilations $\iota'_t(v)=(1-t)v$
$(v\in \Omega',\ t\in [0,1])$ to a family of
holomorphic maps $\iota_t\colon \Omega\to \Omega$ satisfying
the stated properties. (It is possible to choose $\Omega'$ Stein and
with convex fibers: take  
$\Omega'_x = \{v\in \nu_x \colon e^{\phi(x)}|v|^2 <1\}$
where $\phi\colon S\to \R_+$ is a sufficiently fast growing
\spsh\ function and $|v|$ is the Euclidean norm of
the vector $v\in \nu_x\subset \C^N$.)

The above proof extends to the general case:
Let $VT(M)\to M$ denote the {\em vertical tangent bundle of $M$},
consisting of all tangent vector to the (regular)
fibers $M_b=\pi^{-1}(b)$, $b\in B$.
As before we replace $M$ by an open Stein domain containing $S$
as a closed Stein subspace. Then $S_b=S\cap M_b$ is a closed
Stein submanifold of $M_b$ for every $b\in B$,
and $VT(S)$ is a holomorphic vector subbundle of $VT(M)|_S$.
Select finitely many holomorphic sections
$V_1,\ldots, V_N$ of $VT(M)$ that generate the latter bundle
at each point of $M$. (This is possible by a generalization
of Cartan' Theorem A, using an induction on the dimension
of the exceptional set where the sections fail to generate.)
The flow $\phi^j_t$ of $V_j$ is
well defined and holomorphic in a neighborhood of
$M\times \{0\}$ in $M\times \C$, and it remains in the
fibers of $\pi$. We can now complete the proof exactly
as before by  splitting $S\times \C^N=E\oplus \nu$
and restricting the map $F$ (\ref{eqn:flows})
to a suitable Stein neighborhood of the zero section
in $\nu$.
\end{proof}

%
%
%
%
\section{A holomorphic approximation-interpolation theorem}
\label{S4}
In this section we prove the following result which will enable 
us to pass from one stratum to the next one in the proof
of Theorem \ref{CAP:stratified}.

%
%
\begin{proposition}
\label{approximation}
Assume that $X'$ is a closed Stein subvariety of a complex space $X$
and that $K$ is a compact holomorphically convex set in $X$
(as in Theorem \ref{Stein-nbds}) such that $K\cap X'$ is $\cO(X')$-convex.
Let $\pi\colon Z\to X$ be a holomorphic submersion
of a complex space $Z$ onto $X$ (Def.\  \ref{submersion}).
Assume that $U\subset X$ is an open set containing $K$
and $f\colon U \cup X' \to Z|_{U\cup X'}$ is a section whose
restrictions to $U$ and to $X'$ are holomorphic.
Then there exist open Stein neighborhoods $V_j$ of $K\cup X'$ 
in $X$ and holomorphic sections $f_j\colon V_j\to Z|_{V_j}$ $(j=1,2,\ldots)$ such that
$f_j|_{X'}=f|_{X'}$ for all $j\in\N$ and $\lim_{j\to\infty} f_j|_K= f|_K$
(the convergence is uniform on $K$).
\end{proposition}

When $X$ is a Stein manifold, 
$K$ is $\cO(X)$-convex, and $f\colon U\cup X' \to Z$ is a holomorphic
map to a complex manifold $Z$, this is \cite[Theorem 3.1]{FF:EOP}.

\begin{proof}
It suffices to show that for every compact $\cO(X)$-convex set
$L\subset X$ there is a sequence of sections $f_j$ in open neighborhoods
$V_j$ of $K\cup (L\cap X')$ satisfying the conclusion of 
Proposition \ref{approximation}; the result then follows by an
induction over a normal exhaustion of $X$.
Hence we may assume that $X$ is finite dimensional.

Further, by Theorem \ref{Stein-nbds} we can replace $X$ by an open Stein
neighborhood of $K\cup X'$ such that $K$ is $\cO(X)$-convex, and we
can choose the neighborhood $U$ of $K$ to be Stein and relatively compact
in $X$. Since $K$ is $\cO(X)$-convex, it is also $\cO(U)$-convex.

We begin by considering  the case when $Z=X\times \C^p$
and $\pi\colon Z\to X$ is the projection $\pi(x,\zeta)=x$.
We identify sections $X\to Z$ with maps $f\colon X\to \C^p$
to the fiber. By Cartan's extension theorem there is a
holomorphic map $\phi\colon X\to \C^p$ such that
$\phi|_{X'}= f|_{X'}$. There exist finitely many
functions $h_1,\ldots, h_m\in \cO(X)$ that vanish on
the subvariety $X'$ and that generate the ideal sheaf
$\cJ_{X'}$ of $X'$ at every point of $U \Subset X$.
Since $U$ is Stein, Cartan's Theorem B furnishes 
holomorphic maps $g_j\colon U\to \C^p$
$(j=1,\ldots,m)$ such that 
$f=\phi +\sum_{j=1}^m {g_j h_j}$ on $U$.
By the Oka-Weil theorem we can approximate each
$g_j$, uniformly on $K$, by a holomorphic map
$\wt g_j\colon X\to\C^p$. The map 
\begin{equation}
\label{phi-plus}
    \wt f=\phi + \sum{\wt g_j h_j} \colon X\to\C^p
\end{equation}
then approximates $f$ uniformly on $K$, and it agrees with $f$
on $X'$. This gives a sequence satisfying
Proposition \ref{approximation} in this special case.

We now turn to the general case. 
Consider the following subsets of $Z$:
\[
    \wt K=f(K), \ \ \wt X' = f(X'), \ \ \wt U=f(U).
\]
Since $f\colon U\to \wt U$ is biholomorphic,
$\wt U$ is Stein.  Theorem \ref{Stein-nbds}
furnishes an open Stein set $\wt\Omega$ in $Z$
containing $\wt U$ as a closed subvariety.
As $K$ is $\cO(U)$-convex, we infer that $\wt K$ is $\cO(\wt U)$-convex, 
and hence $\cO(\wt \Omega)$-convex. Since $X'$ is Stein
and $K\cap X'$ is $\cO(X')$-convex, $\wt X'$ is a
Stein subvariety of $Z$ and $\wt K\cap \wt X'= f(K\cap X')$
is $\cO(\wt X')$-convex. Theorem \ref{Stein-nbds}
now shows that $\wt K\cup \wt X'$ admits a basis of
Stein neighborhoods in $Z$.

\begin{lemma}
\label{local-embedding}
There exist an open Stein neighborhood $W$
of $\wt K\cup \wt X'$ in $Z$ and a holomorphic embedding
$G\colon W \hookrightarrow X\times\C^N$ for some $N\in \N$
such that for every $x\in X$,
$G$ embeds the fiber $W_x=W\cap \pi^{-1}(x)$
onto a locally closed Stein submanifold
$G(W_x)$ of $\{x\}\times \C^N$.
($W_x$ may be empty for some $x\in X$.)
\end{lemma}

\begin{proof}
Since the fibers $Z_x$ are smooth
and of constant dimension $p$, the vertical tangent bundle
$VT(Z)$, consisting of all tangent vectors to the fibers $Z_x$,
is a holomorphic vector bundle of rank $p$ over $Z$.
For each holomorphic function $g$ on an open subset
of $Z$ we denote by $Vd(g)$ the differential of $g$
in the vertical directions $VT(Z)$; hence $Vd(g)$ 
is a holomorphic section of the vertical cotangent bundle
$VT^*(Z)$.

Let $W_0\subset Z$ be a Stein open neighborhood of $\wt K\cup \wt X'$
furnished by Theorem \ref{Stein-nbds}.
By Cartan's Theorem A there exist finitely many
functions $g_1,\ldots,g_N \in {\cO}(W_0)$ whose vertical
differentials $Vd(g_j)$ span the vertical cotangent
space $VT^*_z(Z)$ at each point $z\in W_0$.
Consider the holomorphic map
\[
     G\colon W_0 \to X\times \C^N,\quad
   G(z)= \bigl( \pi(z), g_1(z),\ldots,g_N(z)\bigr).
\]
Our choice of the $g_j$'s implies that $G$ embeds an
open neighborhood $W\subset W_0$ of $\wt K\cup \wt X'$
biholomorphically onto a locally closed complex subvariety
$G(W)$ of $X\times \C^N$. Clearly we have 
$pr_1\circ\, G=\pi$, where $pr_1\colon X\times \C^N\to X$
is the projection $pr_1(x,\zeta)=x$,
and $G(W_x)$ is a locally closed complex submanifold
(without singularities) of $\{x\}\times \C^N$ for every $x$.
By Theorem \ref{Stein-nbds} we may choose $W$  Stein,
and then each fiber $W_x$ is also Stein.
\end{proof}

We now complete the proof of Proposition \ref{approximation}.
Let $G\colon W\hookrightarrow X\times \C^N$
be a holomorphic embedding furnished by Lemma \ref{local-embedding}.
Proposition \ref{retractions-on-fibers},
applied to the Stein subvariety $S=G(W)$ in $X\times\C^N$
with regular fibers $S_x=G(W_x)\subset \{x\}\times\C^N$,
gives an open Stein neighborhood
$\Omega \subset X\times\C^N$ of $S$
and a fiber preserving holomorphic retraction
$\iota\colon\Omega\to S$ which
retracts each fiber $\Omega_x=\Omega\cap (\{x\}\times\C^N)$
onto the fiber $S_x$.
After shrinking $U$ around $K$ we may assume that
$f(U)\subset W$. Consider the composed section
\[
    G\circ f\colon U\cup X'\to (U\cup X')\times\C^N
    \subset X\times\C^N.
\]
By the special case proved above there is a
sequence of holomorphic sections
$F_j\colon V_j\to V_j \times \C^N$ in open
sets $V_j\supset K\cup X'$ such that
$F_j|_{X'}=G\circ f|_{X'}$ for each $j$, and
$\lim_{j\to \infty} F_j|_K=G\circ f|_K$ uniformly on $K$.
For sufficiently large $j$, and for $V_j\supset K\cup X'$
chosen small enough, we have $F_j(V_j)\subset \Omega$
(the domain of the retraction $\iota$).
The sequence of holomorphic sections
\[
    f_j=G^{-1}\circ \iota\circ F_j\colon V_j\to Z|_{V_j},
    \quad j=1,2,\ldots
\]
then fulfills Proposition \ref{approximation}.
\end{proof}

\begin{remark}
\label{approximation-addition}
Our proof gives the following addition to Proposition \ref{approximation}.
Assume that $f$ satisfies the hypotheses of the proposition, 
and $\phi\colon W\to Z|_W$ is a holomorphic section 
in an open neighborhood $W$ of $X'$ such that $f=\phi$ 
on $X'$ and $f$ and $\phi$ are tangent to order $\ell\in\N$ 
along $K\cap X'$. Then the sequence $\{f_j\}$ in Proposition \ref{approximation} 
can be chosen such that, in addition to the stated properties, 
$f_j$ is tangent to $\phi$ to order $\ell$ along $X'$
for each $j=1,2,\ldots$. 

This follows from the setup (\ref{phi-plus})
by choosing $h_1,\ldots,h_m\in\cO(X)$ that  vanish to order $\ell$
on the subvariety $X'$ and that generate the sheaf $\cJ_{X'}^\ell$ 
(the $\ell$-th power of the ideal sheaf $\cJ_{X'}$) 
at every point of the compact set $K$.
\qed\end{remark}

\begin{remark}
\label{approximation-homotopies}
Proposition \ref{approximation} extends to continuous families of 
sections with a parameter in a compact Hausdorff space. 

As a simple case, let us consider a homotopy of sections
$f_t \colon U\cup X' \to Z$ for $t\in [0,1]$ such that each 
$f_t$ satisfies the conditions in the proposition and the 
homotopy is fixed on $X'$. (Here $U$ is an open set in $X$ containing 
a compact $\cO(X)$-convex set $K$.)
It is then possible to find another homotopy $g_t\colon V\to Z$
$(t\in[0,1])$, consisting of sections that are holomorphic in a neighborhood
$V$ of $K\cup X'$, such that $g_t=f_t=f_0$ on $X'$ and
$g_t$ approximates $f_t$ uniformly on $K$ as close as desired
for every $t\in[0,1]$. In addition, we can choose $g_t$ tangent
to order $\ell\in\N$ along $X'$ to a chosen holomorphic
extension $\phi$ of $f_0|_{X'}$.

This extension is  easy to prove for the trivial fibration $Z=X\times \C^p\to X$.
In the general case a Stein neighborhood of $f_t(K\cup X')$ (in $Z$) is used to get a 
corresponding approximating section $g_t$ for a fixed $t\in [0,1]$, 
and hence a patching problem appears 
when trying to find a family $g_t$ depending continuously on $t\in[0,1]$. 
This problem is solved by the {\em method of successive
patching} explained in \cite[p.\ 139]{FP1}.
\qed\end{remark}

%
%
%
%
\section{Proof of Theorem \ref{CAP:stratified}}
\label{S5}
The proof of Theorem \ref{CAP:stratified} proceeds by a double induction:
The outer one over an exhaustion of $X$ by compact holomorphically
convex sets, and the inner one over the strata in a suitable stratification. 
The main step for the latter inductions is given by the following proposition.

%
%
%
%
\begin{proposition}
\label{special-stratified}
Assume that $X$ is a Stein space,
$M_1\subset M_0$ are closed complex subvarieties of $X$ such that
$S=M_0\bs M_1$ is a complex manifold, and $\pi \colon Z\to X$ is a
holomorphic submersion such that $Z|_S\to S$
is a holomorphic fiber bundle whose fiber enjoys \CAP.
Let $d$ be a complete distance function on $Z$. Given
a pair of compact $\cO(X)$-convex subsets $K\subset L$
of $X$ and a continuous section $f\colon X\to Z$ that is holomorphic
in an open neighborhood of $K_1=K\cup (L\cap M_1)$, there exists
for every $\epsilon>0$ and $\ell \in\N$ 
a homotopy of continuous sections $f_t\colon  X\to Z$
$(t\in[0,1])$ that are holomorphic in a neighborhood of $K_1$, 
with $f_0=f$, satisfying the following:
%
\begin{itemize}
\item[(i)]  $f_t$ agrees with $f_0$ to order $\ell$ along $M_1\cap L$ for each $t\in [0,1]$,
\item[(ii)]  $\sup_{x\in K,\ t\in[0,1]} d\bigl(f_t(x), f_0(x)\bigr) <\epsilon$, and
\item[(iii)] $f_1$ is holomorphic in a neighborhood of
$K_0=K \cup (L\cap M_0)$ in $X$.
\end{itemize}
\end{proposition}

\begin{proof}[Proof of Theorem \ref{CAP:stratified}]
Assume Proposition \ref{special-stratified} for the moment.
Choose a sequence of compact $\cO(X)$-convex sets
\[
    K=K_0\subset K_1\subset K_2\subset \cdots \subset
    \bigcup_{k=0}^\infty K_k =X.
\]
Set $t_k=1-2^{-k}$ and $I_k=[t_k,t_{k+1}]$ for $k=0,1,\ldots$;
thus $\cup_{k=0}^\infty  I_k=[0,1)$.

Let $f=f_0\colon X \to Z$ be a continuous section that 
is holomorphic on a complex subvariety $X'$ of $X$
and in an open neighborhood of $K=K_0$ in $X$. Given $\epsilon >0$,
we need to find a homotopy of sections 
$f_t\colon X\to Z$ $(t\in[0,1])$ satisfying the conclusion (B)
in Theorem  \ref{CAP:stratified}. Essentially the same proof will give the conclusion
(C) if $f_0$ is holomorphic in a neighborhood of $K\cup X'$.

By induction on $k\in \Z_+$ we shall find a sequence of homotopies
of sections $f_t\colon X\to Z$ $(t\in I_k)$ that agree 
at the common endpoint $t_{k+1}$ of the adjacent intervals 
$I_k$, $I_{k+1}$  and such that the following hold:
\begin{itemize}
\item[(i)]  for every $k=0,1,\ldots$ and $t\in I_k$ the section $f_{t}$ 
is holomorphic in an open neighborhood of $K_k$ and it satisfies
\[
	\sup_{x\in K_k} d\bigl(f_t(x),f_{t_k}(x)\bigr) < 2^{-k-1}\epsilon,
\]
\item[(ii)] the homotopy $\{f_t\colon t\in [0,1)\}$ is fixed
on the subvariety $X'$.
\end{itemize}
These properties clearly imply that the limit section
\[
	f_1=\lim_{t\to 1} f_t \colon X\to Z
\] exists and is holomorphic on $X$,
$\sup_{x\in K_0} d\bigl(f_1(x),f_0(x)\bigr)<\epsilon$, and
$f_1|_{X'}= f_0|_{X'}$.  Thus the homotopy $\{f_t\}_{t\in [0,1]}$
satisfies the required properties.

Since all steps in the induction are of the same kind,
we explain how to get the first homotopy for
$t\in I_0= [0,\frac{1}{2}]$. Set $K=K_0$, $L=K_1$.
By the assumption there exists an open set $U\Subset X$
containing $L$ such that $Z|_U$ is a stratified holomorphic
fiber bundle whose strata satisfy CAP. Since $L$ is
$\cO(X)$-convex, there is a Stein domain
$\Omega$ in $X$ with $L\subset \Omega\Subset U$.
Then $\Omega$ is a finite dimensional Stein space and
the restriction $Z|_{\Omega} \to \Omega$ is also a stratified fiber
bundle all of whose fibers enjoy CAP.
Choose a stratification 
\[
	\Omega=X_0 \supset X_1 \supset\cdots\supset  X_m=\emptyset
\]
such that the restriction of 
$\pi\colon Z\to X$ to each stratum $S=X_k\bs X_{k+1}$
is a fiber bundle whose fiber enjoys CAP.
Taking $X'_k= X_k\cup (X'\cap\Omega)$ gives another
stratification
\[
    \Omega =X'_0\supset X'_1\supset\cdots\supset X'_m=\Omega \cap X',
\]
with regular strata $X'_k\bs X'_{k+1}= X_k\bs (X_{k+1}\cup X') \subset X_k\bs X_{k+1}$,
ending with $X'_m=\Omega \cap X'$. 

By Proposition \ref{approximation} we can replace $f_0$ by a section
that is holomorphic in an open neighborhood of $K \cup (L\cap X')$.
When proving (C), $f_0$ is already assumed holomorphic near $X'$.

Let $\{f_t\}_{t\in [0,\frac{1}{2m}]}$ be a homotopy furnished by
Proposition \ref{special-stratified} for the pair of subvarieties
$M_1=X'_m= \Omega\cap X'$ and $M_0= X'_{m-1}$ of $\Omega$,
with $\epsilon$ replaced by $\frac{\epsilon}{2m}$. 
Then $f_{\frac{1}{2m}}$ is holomorphic in a neighborhood
of $K\cup (L\cap X'_{m-1})$, the homotopy is fixed on $X'$, and 
it satisfies
\[ 
	\sup_{x\in K_0}d\bigl(f_t(x),f_0(x)\bigl) < \frac{\epsilon}{2m},\qquad 
  t\in [0,\frac{1}{2m}]. 
\]

Next we apply Proposition \ref{special-stratified}
with the `initial' section $f=f_{\frac{1}{2m}}$ and the pair of subvarieties
$M_1=X'_{m-1}$, $M_0=X'_{m-2}$ to get a homotopy
$\{f_t\}_{t\in [\frac{1}{2m},\frac{2}{2m}]}$ that is fixed on $X'$ 
such that the section $f_{\frac{2}{2m}}=f_{\frac{1}{m}}$ 
is holomorphic in a neighborhood of $K\cup (L\cap X'_{m-2})$, and 
such that
\[ 
	\sup_{x\in K_0} d\bigl(f_t(x),f_{\frac{1}{2m}}(x)\bigl) < \frac{\epsilon}{2m},
  \qquad t\in [\frac{1}{2m},\frac{2}{2m}]. 
\]

Continuing in this way we obtain after $m$ steps
a homotopy $\{f_t\}_{t\in[0,\frac{1}{2}]}$ with the required properties.
In particular, the section $f_{\frac{1}{2}}$ is holomorphic in a
neighborhood of the set $L=K_1$, and its restriction to the subvariety $X'$ 
agrees with $f_0|_{X'}$.
We can extend this homotopy to all of $X$ (without changing it
near $L=K_1$) by using a cut-off function in the parameter.

In the same way we construct homotopies $\{f_t\}_{t\in I_k}$
for all $k=1,2,\ldots$, thereby proving part (B).

When the initial section $f=f_0$ is holomorphic in 
a neighborhood of $K\cup X'$, we can use the improvements
of Proposition \ref{approximation} in 
Remarks \ref{approximation-addition} and \ref{approximation-homotopies} 
to keep the sections $f_t$ in our homotopy holomorphic in a neighborhood 
of $K\cup X'$ and tangent to $f_0$ to order $\ell$ along $X'$.
This gives (C).

Finally we show that the conclusion (D) in Theorem \ref{CAP:stratified}
is a consequence of the interpolation 
statement (B-ii). Suppose that $f_0,f_1\in \Gamma_{\cO}(X,Z)$ are
connected by a homotopy $\{f_t\}_{t\in[0,1]}\subset \Gamma(X,Z)$ 
that is fixed on a closed complex subvariety $X'$ of $X$
(possibly empty). Choose a retraction $\theta$ of $\C$ onto the segment
$[0,1]\subset \C$ and define
$F\colon \wt X= X\times \C\to \wt Z= Z\times\C$ by
\[
	F(x,t)=\bigl(f_{\theta(t)}(x),t\bigr),\qquad x\in X,\ t\in \C.
\]
Then $F$ is a section of the fiber bundle
$\wt\pi \colon \wt Z\to \wt X$ with fiber $Y$ and the projection map
$\wt\pi (z,t)= \bigl(\pi(z),t\bigr)$. 
Since $F$ is holomorphic on the subvariety
\[
	\widetilde X' = (X\times\{0,1\}) \cup (X'\times \C) \subset \wt X,
\]
part (B) of Theorem \ref{CAP:stratified} furnishes a homotopy 
that is fixed on $\widetilde X'$ from $F$ 
to a holomorphic section $\wt F\colon\wt X\to \wt Z$.
The family $\wt f_t\colon X\to Z$, defined by
\[
	\wt F(x,t)=\bigl(\wt f_t(x),t\bigr),\qquad x\in X,\ t\in[0,1],
\]
is then a homotopy of holomorphic sections of $\pi\colon Z\to X$ 
such that $\wt f_0=f_0$, $\wt f_1=f_1$,
and $\wt f_t|_{X'}=f_t|_{X'}$  for all $t\in[0,1]$.

This completes the proof of Theorem \ref{CAP:stratified} provided 
that Proposition \ref{special-stratified} holds.
\end{proof}

\begin{proof}[Proof of Proposition \ref{special-stratified}]
When  $X$ is smooth (without singularities) and $M_0=X$, this is
precisely \cite[Proposition 4.1]{FF:EOP}.
We first recall the proof of this special case
since the general case will be obtained by a modification
explained below. 

We may assume that $L=\{x\in X\colon \rho(x)\le 0\}$,
where $\rho\colon X\to\R$ is a strongly plurisubharmonic exhaustion
function such that $\rho|_K<0$ and $d\rho\ne 0$ on $bL=\{\rho=0\}$.

We recall the geometric setup from
\cite[\S 6.5]{FF:Acta} that was also used in \cite[\S 4]{FF:EOP}.
By the assumption, $f$ is holomorphic in an open set $U\supset K\cup M_1$.
Since the compact set $K_1 = K\cup (M_1 \cap L) \subset U$
is $\cO(X)$-convex, there exists a smooth strongly plurisubharmonic
exhaustion function $\tau \colon X\to \R$ such that $\tau<0$ on $K_1$
and $\tau >0$ on $X\bs U$. By general position we may assume
that $0$ is a regular value of $\tau$ and the hypersurfaces
$\{\rho=0\}=bL$ and $\{\tau=0\}$ intersect transversely along
the real codimension two submanifold $\Sigma = \{\rho=0\}\cap \{\tau=0\}$.
Hence $D_0 :=\{\tau \le 0\} \subset U$ is a strongly pseudoconvex domain
with smooth boundary. For each $s\in[0,1]$ let
\begin{equation}
\label{rhos}
    \rho_s =\tau + s(\rho-\tau)=(1-s)\tau + s\rho,\quad
    D_s=\{\rho_s \le 0\} = \{\tau \le s(\tau-\rho) \}.
\end{equation}
We have $D_0=\{\tau\le 0\}$ and $D_1=\{\rho\le 0\}=L$.

Let $\Omega= \{\rho<0,\ \tau>0\} \subset D_1\bs D_0$
and $\Omega'=\{\rho>0,\ \tau<0\} \subset D_0\bs D_1$ (see Figure \ref{D-s}).
As $s$ increases from $0$ to $1$, $D_s \cap L$  increases
to $D_1=L$ while $D_s\bs L \subset D_0$ decreases to $\emptyset$.
All hypersurfaces $bD_s= \{\rho_s=0\}$ intersect along $\Sigma$. Since
$d\rho_s =(1-s)d\tau+ s d\rho$ and the differentials $d\tau$, $d\rho$
are linearly independent along $\Sigma$, $bD_s$ is smooth near $\Sigma$.
Finally, $bD_s$ is strongly pseudoconvex at every smooth point,
in particular  at every point where $d\rho_s\ne 0$.

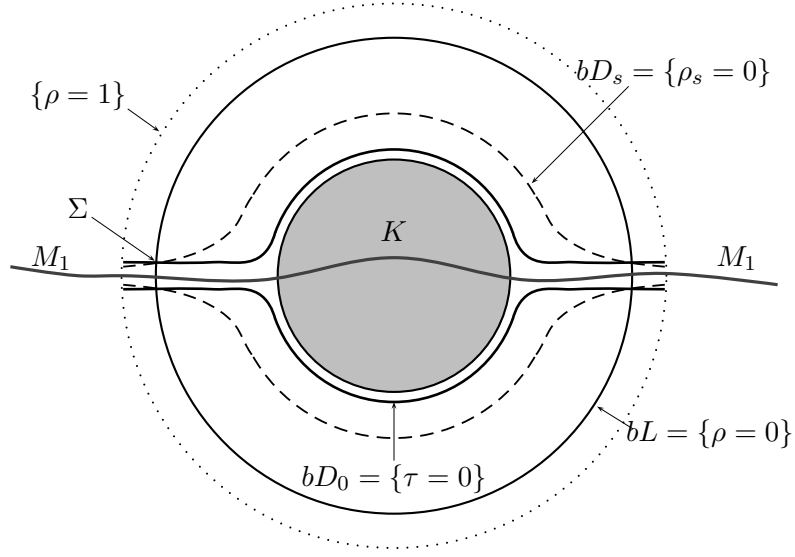
\begin{figure}[ht]
\psset{unit=0.6 cm}
\begin{pspicture}(-8,-6.5)(8,6.5)

%
%
\pscircle[linestyle=dotted](0,0){6.06}                      
\pscircle(0,0){5.3}                                         

%
%
\pscircle[fillstyle=solid,fillcolor=lightgray](0,0){2.6}      

\pscurve[linewidth=1.2pt,linecolor=darkgray](-8.5,0.2)(-7,0)(-6,0)(-3,-0.1)(0,0.4)(3,-0.1)(6,0.05)(8.5,-0.2)
\rput(-7.6,0.4){$M_1$}
\rput(7.6,0.4){$M_1$}

%
%
\psarc[linewidth=1pt](0,0){2.8}{20}{160}
\psarc[linewidth=1pt](0,0){2.8}{200}{340}

%
%
%
\psecurve[linewidth=1pt](3,5)(2.63,0.97)(3,0.4)(4,0.3)(6,0.3)(7,0.3)
\psecurve[linewidth=1pt](3,-5)(2.63,-0.97)(3,-0.4)(4,-0.3)(6,-0.3)(7,-0.3)
\psecurve[linewidth=1pt](-3,5)(-2.63,0.97)(-3,0.4)(-4,0.3)(-6,0.3)(-7,0.3)
\psecurve[linewidth=1pt](-3,-5)(-2.63,-0.97)(-3,-0.4)(-4,-0.3)(-6,-0.3)(-7,-0.3)

%
%
%
\psarc[linestyle=dashed,linewidth=0.7pt](0,0){3.6}{30}{150}
\psarc[linestyle=dashed,linewidth=0.7pt](0,0){3.6}{210}{330}

%
%
%
\psecurve[linestyle=dashed,linewidth=0.7pt](2,3)(3.15,1.75)(3.6,1)(4.4,0.5)(6,0.2)(7,0.2)
\psecurve[linestyle=dashed,linewidth=0.7pt](2,-3)(3.15,-1.75)(3.6,-1)(4.4,-0.5)(6,-0.2)(7,-0.2)
\psecurve[linestyle=dashed,linewidth=0.7pt](-2,3)(-3.15,1.75)(-3.6,1)(-4.4,0.5)(-6,0.2)(-7,0.2)
\psecurve[linestyle=dashed,linewidth=0.7pt](-2,-3)(-3.15,-1.75)(-3.6,-1)(-4.4,-0.5)(-6,-0.2)(-7,-0.2)

\psline[linewidth=0.2pt]{<-}(3.05,2.05)(5.2,4.2)
\rput(6.2,4.5){$bD_s=\{\rho_s =0\}$}

\rput(0,1){$K$}

\rput(-7,1.5){$\Sigma$}
\psline[linewidth=0.2pt]{->}(-6.7,1.3)(-5.3,0.33)

\rput(-7,4){$\{\rho=1 \}$}
\psline[linewidth=0.2pt]{->}(-5.9,3.7)(-5.2,3.2)

\rput(0,-4.5){$bD_0=\{\tau=0\}$}
\psline[linewidth=0.2pt]{->}(0,-4.1)(0,-2.8)

\rput(7,-3.5){$bL=\{\rho=0\}$}
\psline[linewidth=0.2pt]{->}(5.3,-3.4)(4.5,-2.9)

\end{pspicture}
\caption{The sets $D_s$.}
\label{D-s}
\end{figure}

We investigate the singular points of $bD_s=\{\rho_s=0\}$ inside
$\Omega$. (The remaining singular points will be irrelevant.)
The defining equation of $D_s \cap \Omega$ can be written as
$\tau \le s(\tau-\rho)$ and, after dividing by $\tau-\rho>0$, as
\[
    D_s \cap \Omega = \{x\in \Omega \colon h(x)=
    \frac {\tau(x)}{\tau(x)-\rho(x)} \le  s\}.
\]
The critical point equation $dh=0$ is equivalent to
\[
    (\tau -\rho)d\tau - \tau(d\tau -d\rho)= \tau d\rho-\rho d\tau =0.
\]
A generic choice of $\rho$ and $\tau$ insures that there are at most
finitely many solutions $p_1,\ldots,p_m\in \Omega$ and no solution
on $b\Omega$. A calculation shows that at each critical point
the complex Hessians satisfy $(\tau -\rho)^2 H_h = \tau H_\rho - \rho H_\tau$.
Since $\tau>0$ and $-\rho>0$ on $\Omega$, we conclude that $H_h>0$ at such points.
We may assume that distinct critical points of $h$ belong to different level sets.

%
%
We are now ready to prove Proposition \ref{special-stratified}
in the special case.
Let $s_0$ and $s_1$ be two regular values of $h$ on $\Omega$
such that $0\le s_0<s_1 \le 1$ and $h$ has at most one critical point
in $\Omega_{s_0,s_1}= \{x\in\Omega\colon s_0<h(x)< s_1\}$.
Suppose inductively that we have already found a homotopy
$f_t \colon X\to Z$, $t\in [0,s_0]$, satisfying the
conditions in the proposition and such that
$f_{s_0}$ is holomorphic in a neighborhood of
$D_{s_0}$. We wish to deform $f_{s_0}$ to a section
$f_{s_1}$ that is holomorphic in a neighborhood of $D_{s_1}$
by a homotopy consisting of sections that are holomorphic near 
$D_{s_0} \cap D_{s_1}$ and such that the homotopy is fixed 
on $M_1$. We consider two cases.

\medskip
{\em The noncritical case:}
$h$ has no critical values in $\Omega_{s_0,s_1}$.

Recall that a pair $(A,B)$ of compact sets in a complex space $X$ 
is said to be a {\em Cartan pair} if the following hold 
(see \cite[p.\ 695]{FF:CAP}):  
\begin{itemize}
\item[(i)]  $A\cup B$ and $A\cap B$ admit a basis of Stein neighborhoods in $X$, and
\item[(ii)] $\overline {A\backslash B} \cap \overline {B\backslash A}=\emptyset$ 
\ (the separation property). 
\end{itemize}

The main step  is to extend 
a holomorphic section (by appproximation)  to a {\em special convex bump} that we now introduce
(Fig.\ \ref{Fig:specialbump}). The CAP property of the fiber will be used 
(only) at this point of the proof.

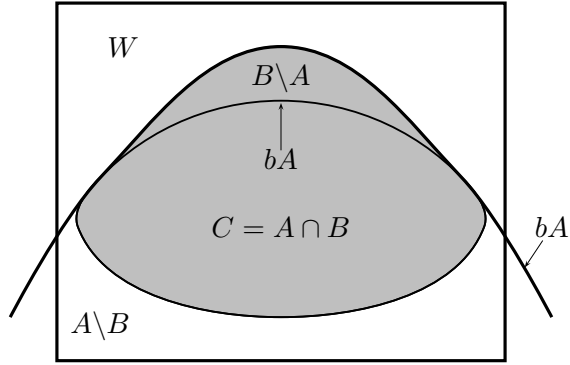
\begin{figure}[ht]
\psset{unit=0.6cm,linewidth=0.7pt}

\begin{pspicture}(-7,0)(7,8)

\pscustom[fillstyle=solid,fillcolor=lightgray]
{
\psecurve(-6,1)(-4,4.1)(0,7)(4,4.1)(6,1)
\psecurve(2,5)(4,4.1)(4.5,3)(0,1)(-4.5,3)(-4,4.1)(-2,5)
}

\psframe[linewidth=1.2pt](-5,0)(5,8)
\pscurve[linewidth=1.2pt](-6,1)(-4,4.1)(0,7)(4,4.1)(6,1)
\psarc[linewidth=0.7pt](0,0.25){5.55}{40}{140}
\psecurve(-2,5)(-4,4.1)(-4.5,3)(0,1)(4.5,3)(4,4.1)(2,5)

\rput(0,3){$C=A\cap B$}
\rput(0,4.5){$bA$}
\psline[linewidth=0.2pt]{->}(0,4.7)(0,5.75)
\rput(0,6.3){$B\bs A$}
\rput(-3.5,7){$W$}
\rput(-4,0.8){$A\bs B$}
\rput(6,3){$bA$}
\psline[linewidth=0.2pt]{->}(5.8,2.7)(5.4,2.1)

\end{pspicture}

\caption{A special convex bump $B$ in a window $W$}
\label{Fig:specialbump}
\end{figure}

%
%
%
%
\begin{definition}
\label{SCP}
Let $X$ be a complex space.
A pair of compact sets $(A,B)$ in $X$ is a {\em special Cartan pair},
and $B$ is a {\em special convex bump on $A$}, if
\begin{itemize}
\item[(i)]    $(A,B)$ is a Cartan pair in $X$, and
\item[(ii)]   there are a compact set $W \subset X_{\rm reg}$
(a {\em window for $B$}),
containing $B$ in its interior, and a holomorphic coordinate map $\phi$
from an open neighborhood of $W$ to $\C^n$ such that $Q = \phi(W)$
is a cube (\ref{cube}), and such that the sets 
\[
	K=\phi(A\cap W),\qquad  K'=\phi((A\cup B)\cap W)
\]
are special convex sets in $Q$ of the form (\ref{special-convex}).
\end{itemize}
\end{definition}

By subdividing $[s_0,s_1]$ into finitely many 
subintervals and replacing $[s_0,s_1]$ by one such subinterval
we can assume that $D_{s_1}$ is obtained by attaching to 
$D_{s_0}\cap D_{s_1}$ finitely many {\em special convex bumps} 
contained in $X\bs M_1$ (see Fig.\ \ref{bump-on-Ds}).
The proof is completely elementary (see e.g.\ \cite[Lemma 12.3]{HL}).

On each bump we apply
\cite[Proposition 3.1]{FF:CAP} to extend the holomorphic homotopy
(by approximation) across the attached bump, without
changing its values on the subvariety $M_1$.
In finitely many steps we accomplish our task.

\begin{figure}[ht]
\psset{unit=0.5 cm}
\begin{pspicture}(-8,-2.7)(8,3)

\psellipse(0,0)(7,2)
\psellipse(0,0)(6,2.7)
\rput(1.8,0.5){$D_{s_0}\cap D_{s_1}$}
\rput(-8.2,1.6){$D_{s_0}$}
\psline[linewidth=0.2pt]{->}(-7.6,1.3)(-6.6,0.7)
\rput(5.3,2.8){$D_{s_1}$}
\psline[linewidth=0.2pt]{<-}(3.5,2.2)(4.7,2.75)

\psecurve[fillstyle=solid,fillcolor=lightgray]
(-0.6,1.9)(0,1.2)(-1,0)(-3,-0.5)(-5,0)(-5.7,0.6)(-5.2,1.3)
(-4,1.95)(-1,1.95)(-0.7,1.9)(0,1.2)(-1,0)  
\rput(-2.8,0.8){$B$}
\psellipse[linestyle=dashed](0,0)(7,2)

\pscurve[linewidth=1.2pt,linecolor=gray](-9,0)(-8,-0.3)
(-7,-0.3)(-5,-0.8)(-3,-1)(0,-0.5)(3,-0.5)(5,-0.2)(7,-0.2)(8,0)
\rput(-7.9,-1.7){$M_1$}
\psline[linewidth=0.2pt]{->}(-8,-1.3)(-8,-0.4)

\end{pspicture}
\caption{A special bump $B$ on $D_{s_0}\cap D_{s_1}$.}
\label{bump-on-Ds}
\end{figure}
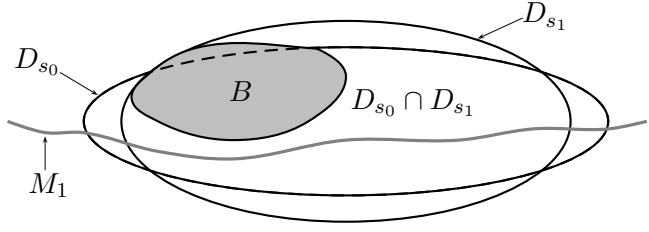

\smallskip
{\em The critical case:}
$h$ has a unique critical point $p\in \Omega_{s_0,s_1}$.
We can extend the holomorphic homotopy across the critical level
$\{h=h(p) \}$ by a reduction to the noncritical case as in
\cite[p.\ 697]{FF:CAP}. This reduction method was developed
in \cite[\S 7.4]{FF:Acta} where the reader can find complete details.

\smallskip
Proposition \ref{special-stratified} is obtained in 
finitely many steps of the above two types.

\smallskip
{\em The general case:}
Now $X$ is a Stein space and $M_1\subset M_0$ are closed
complex (hence Stein) subvarieties of $X$ such that $S=M_0\bs M_1$
is a complex manifold, possibly disconnected. By the assumption, 
we have a section $f\colon X\to Z$ which is holomorphic in 
an open set $U\subset X$ containing $K_1=K\cup (M_1\cap L)$.
As before, we may assume that $L=\{x\in X\colon \rho(x)\le 0\}$,
where $\rho\colon X\to\R$ is a smooth strongly plurisubharmonic
exhaustion function whose restriction to $S=M_0\bs M_1$ 
has no critical points on $bL\cap S=\{\rho=0\}\cap S$. 
Note that $L$ intersects at most finitely many connected
components of $S$, annd hence we can embed a 
relatively compact neighborhood of $L$ in $X$ 
holomorphically into $\C^N$ (see \cite{Bishop}).

Since $K_1$ is  $\cO(X)$-convex, there is a smooth strongly plurisubharmonic
function $\tau \colon X\to \R$ such that $\tau<0$ on $K_1$
and $\tau >0$ on $X\bs U$. Set 
\[
	D_0=\{x\in M_0 \colon \tau(x)\le 0\} \subset U.
\]
By general position we may assume that $0$ is a regular value of 
$\tau|_S$ and the real hypersurfaces $\{\rho=0\} \cap S = bL\cap S$ 
and $\{\tau=0\} \cap S$ in $S$ intersect transversely 
along the submanifold $\Sigma = \{\rho=0\}\cap \{\tau=0\} \cap S$ of $S$.

We define $\rho_s$ as in (\ref{rhos}) and set 
\[
	D_s =\{x\in M_0 \colon \rho_s(x) \le 0\} = 
	\{x\in M_0\colon \tau(x) \le s \bigl(\tau(x)-\rho(x)\bigr) \}.
\]
As $s$ increases from $0$ to $1$, $D_s \cap L$ increases
from $D_0\cap L\subset M_0$ to $D_1=L\cap M_0$. 

Like in the special case considered above, we successively attach 
to the set $A_0=K\cup (D_0\cap L)$ special convex bumps and 
handles, contained in the submanifold $S\cap L$,
thereby reaching the set $K_0=K\cup (L\cap M_0)$ in finitely many steps.
Note that the set $A_0$  is $\cO(X)$-convex, and it contains a 
collar around the set $K_1=K\cup (L\cap M_1)$ in $L\cap M_0$.

We consider a typical step in the {\em noncritical case}.
Assume that $(A,B)$  is a Cartan pair in $X$ 
such that $A$ has been obtained by attaching to $A_0$
finitely many special convex bumps and handles contained in $S\cap L$,
and that $B\subset S\cap L$ is a special convex bump on $A\cap S$
(Def.\ \ref{SCP}). Note that $(A,B)$ is also 
a Cartan pair in the ambient space $\C^N$.

Assume inductively that we have a section $f\colon X\to Z$ that is 
holomorphic in a neighborhood of  $K_1 =K\cup (L\cap M_1)$ in $X$, 
and that is also holomorphic in a relative neighborhood of $A\cap S$ 
in the stratum $S$.
By \cite[Lemma 3.2]{FF:CAP} we `thicken' $f$ over its domain
of holomorphicity to a family of holomorphic sections $F(x,w)$ of $Z$,
depending holomorphically on a parameter $w=(w_1,\ldots,w_p)$ in 
an open neighborhood $O$ of the origin in some Euclidean space $\C^p$,
such that $f=F(\cdotp,0)$, $F(x,w)=f(x)$ for all
$x\in M_1$ and $w \in O$,  and $F$ is submersive
in the $w$-variable for all base points $x$ in a neighborhood of
$A\cap B$ in $S$.

Choose open neighborhoods $O_3 \Subset O_2\Subset O_1$ 
of $0$ in $\C^p$, with $O_1\Subset O$.
By invoking the CAP property of the fiber of the bundle $Z|_S\to S$
we can approximate $F$, uniformly over a neighborhood of $A\cap B$ in $S$
and uniformly with respect to $w\in O_1$, by a family $G$ of holomorphic sections, 
defined in a neighborhood of $B$ in $S$ and depending holomorphically 
on  $w \in O_1$. 

Assuming that $G$ approximates $F$ sufficiently closely,
we find a holomorphic {\em transition map} of the form
\[
	\gamma(x,w)=\bigl(x,\psi(x,w)\bigr),
\]
close to the identity map on $(A\cap B)\times O_2$
(depending on the closeness of $G$ to $F$),
such that we have  
\[ F=G\circ \gamma \quad {\rm on}\ (A\cap B)\times O_2 \]
(see \cite[(3.2)]{FF:CAP}).

We now consider a neighborhood of $L$ in $X$ as a complex subspace of $\C^N$.
Choose a holomorphic retraction $\iota$ from an open
neighborhood $U$ of $A\cap B$ in $\C^N$ onto an open neighborhood
of $A\cap B$ in $S$. (Such retraction exists by Theorem 
\ref{retraction-on-submanifold} since $A\cap B \subset S$
has a Stein neighborhood basis and $S$ is a submanifold of $\C^N$.) 
We choose $U$ small enough such that
$\gamma$ is holomorphic on $\iota(U) \times O_2$. Let 
\[ \wt \gamma(z,w)=\gamma(\iota(z),w),\quad z\in U, \ w \in O''.\]

Choose a Cartan pair $(\wt A,\wt B)$ in $\C^N$ such that
$A\subset {\rm Int}\wt A$, $B\subset {\rm Int}\wt B$, and $\wt A\cap \wt B \subset U$.
(We may take $\wt A$ and $\wt B$ to be smooth strongly pseudoconvex 
domains. See Lemma 4.3 in \cite{FF:Acta}.) 
Assuming as we may that $\wt \gamma$ is sufficiently close to the identity map
(this is insured by approximating $F$ sufficiently well by $G$),
Lemma 2.1 in \cite{FF:CAP} furnishes a splitting
\[
	\wt \gamma = \beta\circ\alpha^{-1},
\]
where $\alpha$ and $\beta$ are biholomorphic maps of
the same type as $\gamma$ and close to the identity
on $\wt A\times O_3$, resp.\ on $\wt B \times O_3$. 
In addition, $\alpha$ can be chosen to match the identity 
map to a given finite order along the intersection of its domain
with the subvariety $M_1$. (For a simple proof of this splitting lemma 
see \cite[Lemma 3.2]{FF:manifolds}.) 

From $F=G\circ\gamma=G\circ\beta\circ\alpha^{-1}$ we infer that  
\[
	F\circ \alpha=G\circ\beta
\]
holds on $(A\cap B)\times O_3$, and this gives a family of holomorphic 
sections $\wt F$ of $Z\to X$ over the set $A\cup B$.
The section $\wt F(\cdotp,0)$ approximates $f$ uniformly on $A$,
it agrees with $f$ to order $\ell$ along $M_1$, and is homotopic to $f$
by a homotopy satisfying the required properties.
The induction may now proceed.

We deal with the critical points of the function 
$h=\frac {\tau}{\tau-\rho}$ in $S \cap L$  
exactly as before by reducing to the noncritical case 
(see the {\em critical case} above).

In finitely many such steps we obtain a homotopy
$\{f_t\}_{t\in[0,1]}$ with the required properties
such that $f_1$ is holomorphic in a neighborhood of
$K_1$ in $X$, and also in a relative neighborhood of 
$L\cap M_0$ in the subvariety $M_0$. 

By Proposition \ref{approximation} there is a holomorphic section
$\wt f_1$ in a neighborhood of $K_0=K\cup (L\cap M_0)$ in $X$ 
such that $\wt f_1$ is as close as desired to $f_1$ on $K$, 
$\wt f_1=f_1$ on $L\cap M_0$, and 
$\wt f_1$ agrees with $f_1$ to order $\ell$ along $L\cap M_1$.
Replace $f_1$ by $\wt f_1$ and adjust the homotopy $\{f_t\}$ accordingly. 
By using a cut-off function in the parameter of the homotopy
we can extend $\{f_t\}$ continuously to all of 
$X$ without changing it near $K_0$ and on $M_1$.

This completes the proof of Proposition \ref{special-stratified}.
\end{proof}

\section{The parametric Oka principle}
\label{parametric}
We have already seen in the proof of
Theorem \ref{CAP:stratified} that the basic Oka principle with 
interpolation also implies the 1-parametric Oka principle,  the latter
meaning that the inclusion
\begin{equation}
\label{inclusion}
    \iota\colon \Gamma_{\cO}(X,Z) \hookrightarrow \Gamma(X,Z)
\end{equation}
of the space of holomorphic sections into the space of continuous sections 
induces a bijection of the path connected components of the two spaces
(see Corollary \ref{CAP-1parametric}).
A stronger form of the Oka principle is to demand that the inclusion
(\ref{inclusion}) is a {\em weak homotopy equivalence}, i.e., 
it induces isomorphisms of all homotopy groups
\begin{equation}
\label{Pik}
	\pi_k(\iota)\colon \pi_k(\Gamma_{\cO}(X,Z)) 
	\stackrel{\approx}{\longrightarrow}  \pi_k(\Gamma(X,Z)),
	\quad k=0,1,2,\ldots
\end{equation}

To formulate the general {\em parametric Oka principle}
we recall  the following notions from \cite[Definition 3]{FP3}.

%
%
\begin{definition}
\label{P-section}
Let $h\colon Z\to X$ be a holomorphic map of complex spaces,
and let $P_0 \subset P$ be topological spaces.
\begin{itemize}
\item[(a)] A {\em $P$-section} of $\pi\colon Z\to X$ is a
continuous map $f\colon X\times  P\to Z$ such that
$f_p=f(\cdotp,p)\colon X\to Z$ is a section of $\pi$ for each fixed
$p\in P$. A $P$-section $f$ is holomorphic on $X$
(resp.\ on a subset $U\subset X$) if $f_p$ is holomorphic on $X$
(resp.\ on $U$) for each fixed $p\in P$.
\item[(b)] A {\em homotopy of $P$-sections} is a $P\times[0,1]$-section,
i.e., a continuous map $H\colon X\times P\times [0,1]\to Z$ such
that $H_t=H(\cdotp,\cdotp,t) \colon X\times P\to Z$ is a
$P$-section for each $t\in [0,1]$. Such homotopy $H$ is holomorphic
if $H_{p,t}=H(\cdotp,p,t) \colon X\to Z$ is holomorphic
for each fixed $(p,t)\in P\times [0,1]$.
\item[(c)] A {\em $(P,P_0)$-section} of $\pi\colon Z\to X$ is a $P$-section
$f\colon X\times  P\to Z$ such that $f_p\colon X\to Z$
is holomorphic on $X$ for each $p\in P_0$.
A $(P,P_0)$-section is holomorphic on $U\subset X$
if $f_p|_U$ is holomorphic for every $p\in P$.
\end{itemize}
A $P$-map $X\to Y$ to a complex space $Y$ is a map
$X\times P\to Y$ (which is the same thing as a $P$-section of the product fibration 
$Z=X\times Y\to X$). Similarly one defines $(P,P_0)$-maps and
their homotopies.
\end{definition}

Given a compact subset $K$ in a complex space $X$ and
a closed complex subvariety $X'$ of $X$, we say that a
$P$-section $f \colon X\times P \to Z$ is holomorphic
on $K\cup X'$ if there is an open neighborhood $U\subset X$
of $K$ such that for every $p\in P$, $f_p= f(\cdotp,p)$ is holomorphic 
in $U$ and $f_p|_{X'}$ is holomorphic on $X'$.
Similar terminology applies to $(P,P_0)$-sections and their homotopies.

%
%
%
%
\begin{definition}
\label{parametric-Oka}
Let $\pi\colon Z\to X$ be a holomorphic map onto a 
complex space $X$, and let $d$ be a distance function on $Z$. 
We say that sections of $\pi$ satisfy the
{\em parametric Oka principle with approximation and
interpolation} {\rm (POPAI)} for a given pair of topological
spaces $P_0\subset P$ if the following holds.

Given a compact $\cO(X)$-convex set $K\subset X$, 
a closed complex subvariety $X'\subset X$, and a $P$-section
$f\colon X\times P\to Z$ that is holomorphic on $K\cap X'$,
there exists for every $\epsilon>0$ a homotopy of $P$-sections
$F\colon X\times P\times [0,1]\to Z$ such that,
setting $f^t=F(\cdotp,\cdotp,t) \colon X\times P\to Z$ $(t\in [0,1])$,
we have 
\begin{itemize}
\item[(i)]   $f^0=f$ is the initial $P$-section,
\item[(ii)]  the $P$-section $f^1$ is holomorphic on $X$,
\item[(iii)] for each $t\in [0,1]$ the $P$-section
$f^t$ is holomorphic on $K\cup X'$, $f^t_p|_{X'} = f^0_p|_{X'}$
for each $p\in P$, 
$\sup_{x\in K,p\in P} d\bigl(F(x,p,t),f(x,p)\bigr) <\epsilon$, and 
\item[(iv)]  if $f$ is a $(P,P_0)$-section then  $F$ can be 
chosen fixed on $P_0$ (i.e., such that $f^t_p$ is independent
of $t\in [0,1]$ for each fixed $p\in P_0$).
\end{itemize}
\end{definition}

Deleting the approximation and/or the interpolation condition 
in POPAI we get POP (the parametric Oka principle), 
POPA (the parametric Oka principle with approximation), 
or POPI (the parametric Oka principle with interpolation).
We also define the 
{\em parametric Oka principle with approximation and jet interpolation}
(POPAJI) by asking that for every $P$-section $f$ which is holomorphic
in an open set $U\supset K\cup X'$ there exist holomorphic $P$-sections 
$\wt f$ that approximate $f$ uniformly on $K$, and that agree with 
$f$ to a given finite order along the subvariety $X'$. 
(The above terminology was introduced by L\'arusson \cite{Larusson3}.)

Note that the validity of POPAJI for all finite polyhedral pairs
$P_0\subset P$ is equivalent to Gromov's 
${\rm Ell}_\infty$ property \cite[\S 3.1]{Gromov}.

The following observation is due to Gromov \cite{Gromov}.

\begin{proposition}
If sections of a holomorphic map  $\pi\colon Z\to X$ satisfy
the parametric Oka principle {\rm (POP)} for all finite polyhedra
$P_0 \subset P$, then the inclusion 
$\iota\colon \Gamma_{\cO}(X,Z) \hookrightarrow \Gamma(X,Z)$
is a weak homotopy equivalence.
\end{proposition}

\begin{proof}
Taking $P=S^k$ (the $k$-sphere) and $P_0=\emptyset$,
POP implies that each continuous map $S^k\to \Gamma(X,Z)$
can be homotopically deformed to a map $S^k\to \Gamma_{\cO}(X,Z)$.
Hence (\ref{Pik}) is surjective.
Similarly, taking $P$ to be the closed real $(k+1)$-ball
$B^{k+1}$ and $P_0=bB^{k+1}=S^k$, we conclude that each map
$S^k\to \Gamma_{\cO}(X,Z)$ that extends to a map
$B^{k+1}\to \Gamma(X,Z)$ also extends to a map
$B^{k+1}\to \Gamma_{\cO}(X,Z)$. Hence
(\ref{Pik}) is injective.
\end{proof}

We now introduce a parametric version of CAP.

\begin{definition}
\label{def:PCAP}
A complex manifold $Y$ enjoys the 
{\em Parametric Convex Approximation Property} {\rm (PCAP)}
for a certain pair of topological spaces $P_0\subset P$ if for every special 
convex pair $(K,Q)$, a map $f\colon Q \times P\to Y$ 
such that $f_p=f(\cdotp,p)\colon Q\to Y$ is holomorphic for every $p\in P_0$,
and is holomorphic on $K$ for every $p\in P$, can be approximated 
uniformly on $K\times P$ by maps $\wt f\colon Q\times P\to Y$
such that $\wt f_p$ is holomorphic on $Q$ for all $p\in P$,
and $\wt f_p=f_p$ for all $p\in P_0$.
\end{definition}

The following is a parametric version of Theorem \ref{CAP:stratified}.

\begin{theorem}
\label{Oka:parametric}
If $\pi\colon Z\to X$ is a stratified holomorphic fiber bundle 
over a Stein space $X$ all of whose fibers satisfy \PCAP\ for a certain
pair of compact Hausdorff spaces $P_0\subset P$, then sections
$X\to Z$ satisfy the parametric Oka principle with approximation and 
(jet) interpolation for the same pair $(P,P_0)$:
\[
	{\rm PCAP}\ \Longrightarrow\ {\rm POPAI},  
	\qquad 	{\rm PCAP}\ \Longrightarrow\ {\rm POPAJI}.
\]
\end{theorem}

\begin{remark}
In previous papers 
\cite{FP1,FP2,FP3,FF:subelliptic,FF:CAP,FF:EOP} POP was only 
considered for pairs of parameter spaces $P_0\subset P$ such that 

(*)  $P$ is a nonempty compact Hausdorff space, and $P_0$ 
is a closed subset of $P$ that is a strong 
deformation neighborhood retract (SDNR) in $P$. 

One can dispense with the SDNR condition by using a generalized
Tietze extension theorem for maps into Hilbert spaces
\cite[Proposition 4.1]{FF:Rothschild} 
\qed \end{remark}

Among the conditions implying PCAP (for all pairs of compact
Hausdorff spaces $P_0\subset P$) are Gromov's 
{\em ellipticity} (the existence of a dominating holomorphic 
spray on $Y$, see \cite{Gromov}),  and {\em subellipticity} 
(the existence of a finite dominating family of holomorphic 
sprays on $Y$, see \cite{FF:subelliptic}). 
Each of these conditions actually implies the parametric Oka principle
for sections of any holomorphic fiber bundle $Z\to X$ with fiber $Y$
over a Stein manifold $X$ (see \cite{FF:subelliptic,FP1,FP2,FP3,Gromov}). 
Like in the basis case considered above,
the proof splits in two parts:
\begin{itemize}
\item subellipticity $\Longrightarrow$ PCAP, and
\item PCAP $\Longrightarrow$ POPAI. 
\end{itemize}

The first of these implications generalizes the 
Oka-Weil approximation theorem (see Theorem 4.1 in \cite[p.\ 135]{FP1}
and Theorem 3.1 in \cite[p.\ 534]{FF:subelliptic}). 
The main part of the proof is the second implication which is indeed
an equivalence (for the converse implication POPAI $\Rightarrow$ PCAP 
it suffices to apply POPAI with $K$ a compact convex set 
in $X=\C^n$ and  $X'=\emptyset$).

\begin{proof}[Sketch of proof of Theorem \ref{Oka:parametric}]
We follow the proof of Theorem \ref{CAP:stratified}, 
using the parametric versions of the technical ingredients
at every step.  There are two types of basic steps 
(see the proof of Proposition \ref{special-stratified}):
\begin{enumerate}
\item extension to a special convex bump (the noncritical case);
\item extension to a handle with totally real core 
(the critical case). 
\end{enumerate}
Step (1) consists of three substeps:
\begin{itemize}
\item[(i)]   holomorphic approximation on a special convex pair,
\item[(ii)]  finding a holomorphic transition map between the old and the new
section,
\item[(iii)] splitting the transition map and gluing the pair of sections. 
\end{itemize}
Each point $p_0\in P$ has an open neighborhood 
$U_{p_0}\subset P$ such that these steps can be performed simultaneously
for all sections $\{f_p\colon p\in U_{p_0}\}$, with a continuous
dependence on the parameter. Indeed, for (i) we use the hypothesis \PCAP,
(ii) follows from the implicit function theorem, 
and (iii) follows from \cite[Lemma 3.2]{FF:manifolds}
that allows continuous dependence on parameters.

Step (2) also involves three substeps:
\begin{itemize}
\item[(a)]   extending the section across the core of the handle,
\item[(b)]  approximation by a section that is holomorphic on a handlebody,
\item[(c)]  passing the critical level by the noncritical
procedure, applied with a different strongly plurisubharmonic function.
\end{itemize}
Substep (a) is automatically fulfilled since we have a
global $P$-section $X\times P\to Z$ at each step.
For Substep (b) we apply the parametric version of \cite[Theorem 3.2]{FF:submersion} 
(see also \cite[Theorem 4.1]{FSlapar1}). 
For Substep (c) we apply Step (1) finitely many times. 

We conclude the proof (globalization with respect to the parameter
$p\in P$) by following the {\em stepwise extension method}
as in the Conclusion of the proof of Theorem 4.2 in \cite{FF:Rothschild}.
\end{proof}

\begin{corollary}
\label{whe}
If $Z\to X$ is a stratified holomorphic fiber bundle over a Stein space $X$
all of whose fibers enjoy {\rm PCAP}, then the inclusion 
\[ 
		\iota\colon \Gamma_{\cO}(X,Z) \hookrightarrow \Gamma(X,Z)
\]
is a weak homotopy equivalence. This holds in particular if
all fibers are elliptic in the sense of Gromov \cite{Gromov},
or subelliptic in the sense of \cite{FF:subelliptic}.
\end{corollary}

\section{Existence of global holomorphic sections}
\label{S7}
We now add a connectivity hypothesis on fibers in
a stratified fiber bundle to obtain existence theorems
for holomorphic sections.

%
%
%
%
\begin{theorem}
\label{stratified-existence}
Assume that $X$ is a reduced finite dimensional Stein space, $\pi \colon Z\to X$ is a
holomorphic submersion, and $X=X_0\supset X_1 \supset\cdots\supset X_m=\emptyset$
is a stratification of $X$ such that for each connected component
$S$ of $X_k\bs X_{k+1}$ the restriction $Z|_S\to S$
is a holomorphic fiber bundle whose fiber $Y_S$ enjoys \CAP\
and $\pi_q(Y_S)=0$ for $q\in\{1,2,\ldots, \dim S-1\}$.
Then there exists a holomorphic section $X\to Z$.

Furthermore, given a closed complex subvariety $X'\subset X$,
a compact $\cO(X)$-convex subset $K \subset X$, an open set $U\supset K$,
and a holomorphic section $f\colon U\cup X' \to Z$, there exists a
holomorphic  section $\wt f\colon X\to Z$ with $\wt f|_{X'}=f|_{X'}$
such that $\wt f$ approximates $f$ as close as desired uniformly on $K$.
\end{theorem}

\begin{proof}
The only place in the proof of
Theorem \ref{CAP:stratified} that requires a
topological condition on the fiber is when crossing
a critical point $p$ of index $k\ge 1$ of a \spsh\ Morse
function $\rho$ on a stratum $S$ (see the {\em critical case} in
the proof of Proposition  \ref{special-stratified}).

To cross the critical level of $\rho$ at $p$
we must be able to extend a given  section,
defined on a sublevel set $\{\rho \le c\}$
for some $c<\rho(p)$ close to the critical level $\rho(p)$,
to a continuous section over a $k$-dimensional totally real disc
$E$ in $S$, attached with its boundary sphere
$bE\approx S^{k-1}$ to the hypersurface $\{\rho=c\}$, 
such that $\{\rho\le c\}\cup E$
is a strong deformation retraction of a sublevel set
$\{\rho\le c'\}$ for some $c'>\rho(p)$.
Such an extension exists if and only if the map $f\colon bE \to Y$
is null-homotopic in $Y$, and this is certainly the case if
the group $\pi_{k-1}(Y)$ vanishes.

Since $\rho|_S$ is strongly plurisubharmonic,
we have $k\le \dim S$.  Our topological assumption on the fibers of $Z\to X$
therefore insures the existence of a continuous extension
of a section at each critical point on every stratum,
and hence the proof of Theorem \ref{CAP:stratified} applies.
\end{proof}

\begin{corollary}
Assume that $X$ is a Stein manifold of dimension $n$,
$K$ is a compact $\cO(X)$-convex subset of $X$, $X'$
is a closed complex subvariety of $X$, $U$ is an open set
in $X$ containing $K$,
and $f\colon U\cup X' \to Y$ is a holomorphic map.
If $Y$ satisfies \CAP\ and $\pi_q(Y)=0$ for $q=0,1,\ldots,n-1$
then there exists a holomorphic map
$\wt f\colon X\to Y$ such that $\wt f|_{X'}=f|_{X'}$
and $\wt f$ approximates $f$ uniformly on $K$ as close as desired.

The conclusion holds if $Y=\C^N\bs A$, where $A$ is a closed
algebraic subvariety of $\C^N$ of codimension
$q=N-\dim A \ge \max\{2,\frac{n+1}{2}\}$.
\end{corollary}

\begin{proof}
The first conclusion is a special case of Theorem \ref{stratified-existence}.

For the second part, recall that $\C^N\bs A$ enjoys CAP if $q\ge 2$
(see \cite[Corollary 1.3]{FF:CAP}), and $\pi_k(\C^N\bs A)=0$
when $0\le k\le 2q-2$ (this follows from a general position 
argument). 
\end{proof}

\section{Submersions with stratified sprays}
Assume that $\pi\colon Z\to X$ is a holomorphic submersion 
of (reduced, finite dimensional)
complex spaces (Def.\ \ref{submersion}). For any point 
$x\in X$ we set $Z_x=\pi^{-1}(x)$. For $z\in Z_x$ 
we denote by $VT_z Z= T_z Z_x$ the tangent space to the 
fiber $Z_x$, also called the {\em vertical tangent space}
of $Z$ at $z$ (with respect to $\pi$).

\begin{definition}
\label{fiber-spray}
({\rm Gromov \cite{Gromov}})
Assume that $\pi\colon Z\to X$ is a holomorphic submersion 
of complex spaces. A {\em fiber dominating spray} associated to 
this submersion is a triple $(E,p,s)$, where $p\colon E\to Z$ is a holomorphic 
fiber bundle and $s\colon E\to Z$ is a holomorphic map satisfying the 
following properties for every $z\in Z$:
\begin{itemize}
\item[(i)] $s(0_z)=z$ (here $0_z\in E_z$ is the zero element of 
the fiber of $E$ over $z$),
\item[(ii)] $\pi\circ s=\pi\circ p$ (i.e., $s(E_z)\subset Z_{\pi(z)}$), and
\item[(iii)] $(ds)_{0_z}(E_z)= VT_z Z$.
\end{itemize}
\end{definition}

Condition (iii), which means that the map $s\colon E_z\to Z_{\pi(z)}$ 
is a submersion at $0_z\in E_z$, 
is called the {\em domination property} of the spray.

\begin{definition}
\label{stratified-spray}
A holomorphic submersion $\pi\colon Z\to X$ admits 
{\em stratified sprays} (or is a {\em stratified elliptic submersion})
if there exists a stratification 
$X=X_0\supset X_1\supset\cdots\supset  X_m=\emptyset$
such that each point $x$ in any stratum $S\subset X_k\bs X_{k+1}$ 
has an open neighborhood $U$ in $S$ 
such that the restricted submersion $\pi\colon Z|_U\to U$ admits a 
fiber dominating spray (Def.\ \ref{fiber-spray}). 
\end{definition}

A special case of the following result was given in \cite[p.\ 66]{FP3}.

%
%
%
%
\begin{theorem}
\label{stratified-sprays}
Let $\pi\colon Z\to X$ be a holomorphic submersion of reduced 
finite dimensional complex spaces. If $X$ is Stein and if
the submersion admits stratified sprays (Def.\ \ref{stratified-spray})
then sections $X\to Z$ satisfy the parametric Oka principle 
with approximation and (jet) interpolation 
for any pair of compact Hausdorff spaces $P_0\subset P$
(see Def.\ \ref{parametric-Oka}). 
The conclusion of Corollary \ref{whe} holds as well.
\end{theorem}

The validity of Theorem \ref{stratified-sprays}  
was indicated by Gromov \cite{Gromov}.
In \cite{FP2,FP3} a complete proof was given for the
case when the base $X$ is a Stein manifold and there is only
one stratum. The  methods of this paper 
allow us to complete the proof also for Stein spaces with 
singularities, and for submersions with stratified sprays,
along the lines indicated in \cite[\S 7]{FP3}. 

\begin{proof}
We follow the proof of Theorem \ref{CAP:stratified}
(see \S 5 above). Lets us consider the basic case without parameteres.
At a typical inductive step we have a section 
$f\colon X\to Z$ that is holomorphic on a neighborhood of a 
compact $\cO(X)$-convex subset $K$ of $X$, and whose restriction 
to a closed complex subvariety $X'$ of $X$ is holomorphic on $X'$.
Given a larger compact $\cO(X)$-convex subset $L$ in $X$ containing $K$,
the goal is to homotopically deform $f$ to a section that is holomorphic 
in a neighborhood of $L$ such that the homotopy is fixed on $X'$
and such that all sections in the homotopy 
are holomorphic near $K$ and are uniformly close to $f$ on $K$.
The proof is then completed by induction over a sequence
of compacta exhausting $X$.

Applying Proposition \ref{approximation} we first modify $f$ to make it holomorphic 
in an open neighborhood of $K\cup X'$. 
Let $X=X_0\supset X_1\supset\cdots\supset  X_m=\emptyset$ be a stratification  
satisfying the hypotheses of the theorem. Replacing $X_k$ by 
$X'_k = X_k\cup X'$ we obtain another stratification 
\[
    X = X'_0\supset X'_1\supset\cdots\supset X'_m= X'
\]
with regular strata $X'_k\bs X'_{k+1}= X_k\bs (X_{k+1}\cup X') \subset X_k\bs X_{k+1}$
such that the dominating spray condition still holds over the strata.  

Assume inductively that $f$ is already holomorphic in a neighborhood of 
$K\cup (X'_{k+1} \cap L)$ (this holds when $k+1=m$).  The goal is to modify
$f$ to a section that is holomorphic in a neighborhood of $K\cup (X'_k\cap L)$
such that the deformation is fixed on $X'_{k+1}$ and is small on $K$;
the proof is then completed by a finite downward induction on $k$
by reaching $k=0$. This is accomplished by Proposition 6.1 in \cite{FP3}
as follows: 

We choose a {\em Cartan string} $(A_0,A_1,\ldots,A_n)$ in 
the subvariety $X'_k$ such that $A_0$ is a (small) neighborhood of 
$K\cup (X'_{k+1}\cap L)$ on which $f$ is holomorphic, while the remaining 
sets $A_1,\ldots,A_n$ are contained in the strata $S=X'_{k}\bs X'_{k+1}$ 
(hence in the regular locus of $X'_k$), and they are chosen 
sufficiently small such that the restricted submersion $Z|_S\to S$ 
admits a fiber dominating spray on a neighborhood of each of them. 
In addition, we require that $\bigcup_{j=0}^n A_j = L\cap X'_k$. 
(For the definition and the construction of such Cartan strings 
see Gromov \cite[4.2.D']{Gromov} and \cite{FP2}, 
Def.\ 4.1 and Theorem 4.6.)
Now we obtain a desired modification of $f$ by applying 
Proposition 5.1 in \cite{FP2} with the given Cartan string.

The proof in the parametric case follows the same lines
(see \S \ref{parametric} above and \cite{FP2,FP3}).
\end{proof}

\begin{remark}
\label{finalremark}
Two further improvements of Theorem \ref{stratified-sprays} are possible. 

As in Theorem \ref{CAP:stratified},
it suffices that the assumptions of Theorem \ref{stratified-sprays} hold 
over an exhausting sequence of open relatively compact subsets of $X$.
This allows $X$ to be infinite dimensional. 

Secondly, we can weaken the fiber dominating condition by asking that each
point $x$ in a statum $S$ (in a given stratification) 
has an open neighborhood $U\subset S$ such that
the restricted submersion $Z|_U\to U$ admits a finite 
collection of holomorphic sprays $(E_j,p_j,s_j)$, $j=1,\ldots,k$,
satifying the following domination condition 
at every point $z\in Z|_U=\pi^{-1}(U)$:
\[
	(ds_1)_{0_z}(E_{1,z}) + (ds_2)_{0_z}(E_{2,z}) + \ldots +
	 (ds_k)_{0_z}(E_{k,z}) = VT_z Z.
\]
We refer to \cite{FF:subelliptic} for further details. 
\end{remark}

\medskip
{\em Acknowledgement.}
I wish to thank Jasna Prezelj for stimulating discussions
on this subject over many years, and Frank Kutzschbauch for 
telling me about the solution to the holomorphic Vaserstein problem 
which prompted me to include \S 8 to the paper.
Finally, I thank the referee for his helful remarks 
that lead me to include a discussion of the parametric case.

\bibliographystyle{amsplain}

\end{document}